\documentclass{amsart}

\usepackage[psamsfonts]{amssymb}

\theoremstyle{plain}
\newtheorem{theorem}{Theorem}[section]
\newtheorem{lemma}[theorem]{Lemma}
\newtheorem{proposition}[theorem]{Proposition}
\newtheorem{cor}[theorem]{Corollary}

\theoremstyle{definition}
\newtheorem{defn}{Definition}[section]

\theoremstyle{remark}

\title[A Rodrigues formula in two variables]{A matrix Rodrigues formula for classical orthogonal polynomials in
two variables}

\author[M. \'{A}lvarez de Morales]{Mar\'{\i}a \'{A}lvarez de Morales}
\address[M. \'{A}lvarez de Morales]{Departamento de Ma\-te\-m\'a\-ti\-ca Apli\-ca\-da,
Uni\-ver\-si\-dad de Gra\-na\-da,
 Gra\-na\-da,
Spain} \email{alvarezd@ugr.es}

\author[L. Fern\'andez]{Lidia Fern\'andez}
\address[L. Fern\'andez]{Departamento de Ma\-te\-m\'a\-ti\-ca Apli\-ca\-da,
Uni\-ver\-si\-dad de Gra\-na\-da,
 Gra\-na\-da,
Spain} \email{lidiafr@ugr.es}

\author[T. E. P\'erez]{Teresa E. P\'erez}
\address[T. E. P\'erez]{Departamento de Matem\'atica Aplicada,
  and Instituto Carlos I de F\'{\i}sica Te\'orica y Computacional,
  Universidad de Granada,
Granada,
 Spain}
\email{tperez@ugr.es}

\author[M. A. Pi\~{n}ar]{Miguel ~~ A. ~~ Pi\~{n}ar}
\address[M. A. Pi\~{n}ar]{Departamento de Matem\'atica Aplicada,
  and Instituto Carlos I de F\'{\i}sica Te\'orica y Computacional,
   Universidad de Granada,
Granada,
 Spain}
\email{mpinar@ugr.es}

\thanks{Partially supported by Ministerio de Ciencia y Tecnolog\'{\i}a
        (MCYT) of Spain and by the European Regional Development Fund
        (ERDF) through the grant MTM 2005--08648--C02--02, and Junta de
        Andaluc\'{\i}a, Grupo de Investigaci\'on FQM 0229.}

\begin{document}



\begin{abstract}
Classical orthogonal polynomials in one variable can be
characterized as the only orthogonal polynomials satisfying a
Rodrigues formula. In this paper, using the second kind Kronecker
power of a matrix, a Rodrigues formula is introduced for classical
orthogonal polynomials in two variables.
\end{abstract}

\maketitle

\section{Introduction}

One of the most important characterizations for classical
orthogonal polynomials in one variable (Hermite, Laguerre, Jacobi
and Bessel) is the so--called {\it Rodrigues formula} (see, for
instance, \cite{Chi}).

Using this kind of formula we can write the $n$--th classical
orthogonal polynomial in terms of a $n$--th order derivative. In
fact, if we denote by $\{P_n\}_n$ a classical family of orthogonal
polynomials in one variable, then
\begin{equation}\label{1rod}
P_n(x) = \frac{k_n}{\omega(x)}\frac{d^n}{dx^n}
\left(\phi(x)^n\,\omega(x)\right), \quad n=0, 1, 2, \ldots,
\end{equation}
where $k_n$ is a constant, $\phi(x)$ is a polynomial of degree less
than or equal to 2, independent of $n$, and $\omega(x)$ is an
integrable function in a appropriate support set.

If $\deg \phi=0$, Hermite polynomials appear, up to a linear change
in the variable. If $\deg\phi=1$, Laguerre polynomials are obtained,
and if $\deg\phi=2$, we can deduce two families of polynomials,
Jacobi polynomials when $\phi(x)$ has two simple roots, and Bessel
polynomials when $\phi(x)$ has a double root.

Formula (\ref{1rod}) is called {\it Rodrigues formula}, honoring
B. O. Rodrigues who established the formula in 1814 for Legendre
polynomials.

\bigskip

Orthogonal polynomials in two variables which are solutions of
partial differential equations were systematically studied by H. L.
Krall and I. M. Sheffer (\cite{KS}), in 1967. They defined {\it
classical} orthogonal polynomials in two variables as the sequences
of orthogonal polynomials $\{P_{h,k}\}_{h,k\ge0}$ such that every
polynomial $P_{h,k}$, with $h+k=n$, satisfies the second order PDE
\begin{equation}\label{KSorig}
L[w]\equiv a \, w_{x x}
 + 2 \, b \, w_{x y}
 + c \, w_{y y}
 + d \, w_x + e \, w_y = \lambda_n w,
\end{equation}
where $a(x,y) = a x^2 + d_1 x + e_1 y +f_1$; $b(x,y) = a x y + d_2 x
+ e_2 y +f_2$; $c(x,y) = a y^2 + d_3 x + e_3 y + f_3$; $d(x,y) = g x
+ h_1$; $e(x,y) = g y + h_2$, and $\lambda_n = a n (n-1) + g n$.

The special shape of the polynomials involved in the above
equation is a direct consequence of the fact that every orthogonal
polynomial of total degree $n$ must satisfy the same PDE. Krall
and Sheffer showed that, up to a linear change in the variables,
there are nine different sets of orthogonal polynomials satisfying
such type of PDE.

The first reference to a Rodrigues formula for classical orthogonal
polynomials in two variables appears in the classical text by P.
Appell and J. Kamp\'{e} de F\'{e}riet (\cite{AK}). Later, P. K. Suetin
(\cite{Su}), and Y. J. Kim, K. H. Kwon and J. K. Lee (\cite{KKL2})
consider an analogue of the Rodrigues formula for Krall and Sheffer
classical orthogonal polynomials in two variables. In fact, for $n$
a positive integer, they define
 \begin{equation}\label{rodorig} P_{n-i,i}(x,y) =
\frac{1}{~\omega~} \,\partial^{n-i}_x\,\partial_y^{i}(p^{n-i}\,
q^i\,\omega), \end{equation}
 where $w(x,y)$ is a weight function over a simply connected domain, and a symmetry factor
 of $L$, the linear differential operator defined in
(\ref{KSorig}), and $p(x,y)$, $q(x,y)$ are polynomials related
with the polynomial coefficients in (\ref{KSorig}). Then, under
some additional hypothesis, (\ref{rodorig}) defines an algebraic
polynomial in two variables orthogonal to all polynomials of lower
degree (see example 3, in Section \ref{mainsec}).

The above Rodrigues formula runs only for classical orthogonal
polynomials associated with a positive definite moment functional,
since it needs a weight function. Nevertheless, H. L. Krall and I.
M. Sheffer founded classical orthogonal polynomials in two
variables associated with a non positive definite moment
functional which has a symmetry factor (see L. L. Littlejohn
\cite{Li}), but not a Rodrigues formula like (\ref{rodorig})
(\cite{KKL2}).

On the other hand, tensor product of two classical orthogonal
polynomials in one variable, defined by
$$
P_{h,k}(x,y)=R_{h}(x) S_{k}(y), \quad h, k \ge 0,
$$
where $\{R_h\}_{h\ge0}$ and $\{S_k\}_{k\ge0}$ are Hermite, Laguerre,
Jacobi or Bessel polynomials, satisfies a Rodrigues formula as
(\ref{rodorig}). In fact, $P_{h,k}(x,y)$ can be written as a product
of the respective Rodrigues formulas
$$P_{h,k}(x,y) = \frac{1}{~\omega_1~}
\partial_x^{h} \,(\phi_1^{h}\,\omega_1) \,\frac{1}{~\omega_2~} \, \partial_y^{k}\,
(\phi_2^{k}\,\omega_2) , \quad h, k \ge 0.$$

However, tensor products of classical orthogonal polynomials in one
variable are not classical according to the Krall and Sheffer
definition since they do not satisfy equation (\ref{KSorig}), except
for Hermite and Laguerre polynomials.

\bigskip

Recently, the authors (see \cite{FPP1,FPP2,FPP3,FPP4}) extended the
concept of classical orthogonal polynomials in two variables to a
wider framework, which, of course, includes the Krall and Sheffer
definition and tensor products of classical orthogonal polynomials
in one variable.

The vector representation for orthogonal polynomials introduced in
\cite{Ko1,Ko2}, and developed in \cite{Xu2} is the key to
introduce the concept of {\it classical} orthogonal polynomials in
two variables. Let $\{\mathbb{P}_n\}_n$ denote a {\it weak
orthogonal polynomial sequence} (see Section 2), it will be called
{\it classical} (in an extended sense) if there exist non singular
matrices $\Lambda_n\in{\mathcal M}_{n+1}(\mathbb{R})$, such that,
\begin{equation}\label{defintro}
L[\mathbb{P}^t_n]\equiv\hbox{\rm div~}(\Phi\nabla
\mathbb{P}^t_n)+\tilde{\Psi}^t\nabla\mathbb{P}^t_n = \mathbb{P}^t_n
\Lambda_n, \end{equation}
 where
$$\Phi = \left(%
\begin{array}{cc}
  a & b \\
  b & c \\
\end{array}%
\right), \qquad \tilde{\Psi} = \left(%
\begin{array}{c}
  d - a_x - b_y \\
  e - b_x - c_y \\
\end{array}%
\right),$$ and $a, b, c$ are polynomials in two variables of total
degree less than or equal to 2, and $d, e$ are polynomials in two
variables of total degree less than or equal to 1, and $\hbox{\rm
div~}$ and $\nabla$ denote the usual divergence and gradient
operators in two variables. Observe that the left hand side of
(\ref{defintro}) generalizes the left hand side of the Krall and
Sheffer PDE (\ref{KSorig}), without any restrictions on the
polynomial coefficients. Moreover, this new definition also
include the tensor product of classical orthogonal polynomials in
one variable. In the Krall and Sheffer case the matrices
$\Lambda_n$ are scalar matrices, and in the tensor product case,
they are diagonal non--singular matrices.

In this paper, we will obtain a matrix Rodrigues type formula for
classical orthogonal polynomials in {\it extended sense}. Denoting
by $\Phi^{\{n\}}$ the second kind Kronecker power of the matrix
$\Phi$ (see Bellman \cite{Be}), we will show, under some
hypothesis, that the expression,
\begin{equation}\label{Rodnues}
\mathbb{Q}_n^t = \frac{1}{~\omega~}\, {\hbox{\rm div}}^{\{n\}}
(\Phi^{\{n\}}\,\omega),\qquad n\ge 0, \end{equation}
 provides a
classical WOPS, where $\omega(x,y)$ is a symmetry factor of the PDE
(\ref{defintro}), and ${\hbox{\rm div}}^{\{n\}}$ is a $n$--th order
differential operator.

This formula generalizes in a natural way the Rodrigues formula
proved in \cite{Su}, and \cite{KKL2}, and the Rodrigues formula
for tensor product of classical orthogonal polynomials in one
variable.

Moreover, using our results, we will deduce a matrix Rodrigues
formula for classical orthogonal polynomials associated with a non
positive definite moment functional whose PDE has a symmetry
factor (see example 6, in Section \ref{mainsec}).

The structure of the paper is as follows. In Section 2 we collect
the necessary basic tools. Section 3 and 4, are devoted to introduce
classical orthogonal polynomials in two variables and symmetry
factors associated with the partial differential equation
(\ref{defintro}). The matrix Rodrigues formula (\ref{Rodnues}), as
well as some examples are studied in Section 6, and finally, the
proof of the main result is given in the last section.

\section{Orthogonal polynomials in two variables}

First, we introduce some notations. Let $\mathcal{P}$ denote the
linear space of real polynomials in two variables, and ${\mathcal
P}_n$ the subspace of polynomials of total degree not greater than
$n$.

Let ${\mathcal M}_{h\times k}(\mathbb{R})$ and ${\mathcal
M}_{h\times k}(\mathcal{P})$ denote the linear spaces of $h\times
k$ real and polynomial matrices, respectively. When $h=k$, the
second index will be omitted.

Let $A$ be a matrix, we denote by $A^t$ its transpose, and by
$\det(A)$ its determinant. As usual, we say that $A$ is
non--singular if $\det(A)\neq 0$. Furthermore, we introduce $I_h$
as the identity matrix of dimension $h$.

Moreover, we define the \emph{degree of} a matrix of polynomials
$A\in\mathcal{M}_{h\times k}(\mathcal{P})$, as
$$\deg A = \max\{\deg a_{i,j}(x,y), 1\le
i\le h, 1\le j\le k\}\ge 0,$$
where $a_{i,j}(x,y)$ denotes the $(i,j)$--entry of $A$.

\bigskip

Before discussing our approach, we briefly give some general
properties and tools about bivariate orthogonal polynomials. For an
exhaustive description of this and another related subjects see, for
instance, \cite{DX,KKL1,KKL2,Koor,Ko1,Ko2,Su,Xu2}.

Let $\{\mu_{h,k}\}_{h,k\ge0}$ be a double indexed sequence of real
numbers, and let $u: {\mathcal P}\rightarrow\mathbb{R}$ be a
functional defined by means of the moments $\mu_{h,k} = \langle
u,x^h\, y^k\rangle$, $h, k = 0, 1, 2,\ldots$, and extended by
linearity. Then, we will say that $u$ is a \emph{moment functional}.

Some elementary properties about moment functionals acting over
polynomial matrices $A \in {\mathcal M}_{h\times k}$ and $B\in
    {\mathcal M}_{h\times l}$ are given by (see
\cite{DX,KKL1,KKL2,Xu2}),
\begin{enumerate}
\item $\langle u, A \rangle =\left(\langle u, a_{i,j}\rangle
\right)_{i,j=1}^{h,k}\in {\mathcal M}_{h\times k}(\mathbb{R}),$\,
where $A =\left(a_{i,j}\right)_{i,j=1}^{h,k}$,
\item $\langle A u, B \rangle = \langle u, A^t~B\rangle.$
\end{enumerate}

We say that a polynomial $p(x,y) \in {\mathcal P}_n$ is {\it
orthogonal} with respect to $u$ if
$$\langle u, p\, q\rangle = 0, \quad \forall q\in {\mathcal
P}, \quad \deg q < \deg p.$$
Then, we can define
$${\mathcal V}_n = \{ p\in{\mathcal P}_n/ \langle u, p\, q\rangle = 0, \forall q\in {\mathcal
P}_{n-1}\}.$$ A moment functional $u$ is called {\it quasi
definite} if $\dim {\mathcal V}_n = n+1$.

\begin{defn}[\cite{KKL1}] A {\emph polynomial system (PS)} is a vector sequence
$\{\mathbb{P}_n\}_{n\geq 0}$ such that
$$\mathbb{P}_n = (P_{n,0},P_{n-1,1}, \ldots, P_{0,n})^t \in {\mathcal M}_{{(n+1)}\times 1}({\mathcal P}_n),
$$
where $\{P_{n,0},P_{n-1,1}, \ldots, P_{0,n}\}$ are polynomials of
total degree $n$ independent modulus $\mathcal{P}_{n-1}$.
\end{defn}

Observe that a PS is a sequence of vectors whose dimension and
total degree are increasing: $\mathbb{P}_0$ is a constant,
$\mathbb{P}_1$ is a column vector of dimension $2$ of bivariate
polynomials of total degree 1, $\mathbb{P}_2$ is a column vector
of dimension $3$ whose elements are bivariate polynomials of total
degree 2, and so on.

\begin{defn}[\cite{KKL1}]
We will say that a PS  $\{\mathbb{P}_n\}_{n\geq 0}$ is a
\emph{weak orthogonal polynomial system (WOPS)} with respect to a moment functional
$u$ if
\begin{eqnarray*}
   \langle u,{\mathbb{P}}_{n}{\mathbb{P}}_{m}^t\rangle &=& 0, \quad n\neq m,\\
   \langle u,{\mathbb{P}}_{n}{\mathbb{P}}_{n}^t\rangle &=& H_n, \quad n = 0, 1, 2, \ldots
\end{eqnarray*}
where $H_n\in {\mathcal M}_{n+1}(\mathbb{R})$ is a non--singular
matrix.
\end{defn}
In the particular case where $H_n$ is a diagonal matrix, we will
say that the WOPS $\{\mathbb{P}_n\}_{n\ge0}$ is an
\emph{orthogonal polynomial system (OPS)}. Moreover, if $H_n =
I_{n+1}$, we call $\{\mathbb{P}_n\}_{n\geq 0}$ an {\it orthonormal
polynomial system}. A moment functional $u$ is \emph{quasi
definite} if and only if there exists a WOPS with respect to $u$
(\cite{KKL1}).

In addition, a WOPS is called a \emph{monic} WOPS if every
polynomial contains only one monic term of higher degree, that is,
$$P_{h,k}(x,y)= x^h \, y^k + R(x,y), \quad h + k = n,$$
where $R(x,y)\in{\mathcal P}_{n-1}$. And finally, we have that for
a quasi definite moment functional $u$, there exists a unique
monic WOPS associated with $u$.

In this paper, we will need some differentiation tools. In fact, we
will use the {\it gradient operator} $\nabla$, and the {\it
divergence operator} $\hbox{\rm div}$, defined as usual. The
extension of this operators for matrices is introduced in
\cite{FPP1,FPP2,FPP3,FPP4}. Let $A, B_0, B_1 \in {\mathcal
M}_{h\times k}({\mathcal P})$ be polynomial matrices. We define
$$
\nabla A = \begin{pmatrix}
                \partial_x A \cr
                \partial_y A\end{pmatrix} \in {\mathcal M}_{2h\times k}({\mathcal
            P}),\quad
\hbox{\rm div~} \begin{pmatrix}B_0 \cr B_1\end{pmatrix} =
\partial_x B_0 +
\partial_y B_1 \in {\mathcal M}_{h\times k}({\mathcal P}).
$$
We extend these definitions for $n\ge 1$. In fact, writing
$\nabla^{\{1\}}= \nabla$, and $\hbox{\rm div}^{\{1\}} = \hbox{\rm
div}$, if we denote $\mathcal{D}_i^n = \binom{n}{i}
\partial_x^{n-i}\,\partial_y^i$, $i=0,1,\ldots,n$, we can
introduce the differential operators $\nabla^{\{n\}}$ and
$\hbox{\rm div}^{\{n\}}$ by means of
\begin{eqnarray}
&~& \nabla^{\{n\}} A = (\mathcal{D}_0^n \, A, \mathcal{D}_1^n \,
A, \cdots, \mathcal{D}_n^n \, A)^t \in {\mathcal
M}_{((n+1)h)\times k}({\mathcal
            P}),\label{nablan}\\
&~&\hbox{\rm div}^{\{n\}} (B_0, B_1, \cdots, B_n)^t = \sum_{i=0}^n
\, \mathcal{D}_i^n\, B_i \in {\mathcal M}_{h\times k}({\mathcal
P})\label{divn},
\end{eqnarray}
where $A, B_0, B_1, \ldots, B_n \in {\mathcal M}_{h\times
k}({\mathcal P})$ are polynomial matrices. In addition, we establish
$\nabla^{\{0\}} A = A$, and ${\hbox{\rm div}}^{\{0\}} A = A$.

The previous definitions can be translated to the linear space of
moment functionals using duality. For $n\ge 0$, we define the {\it
n--th distributional gradient operator} and the {\it n--th
distributional divergence operator} acting over moment functionals
in the following way
\begin{eqnarray*}
&~& \langle\nabla^{\{n\}} u , \begin{pmatrix}p_0\cr p_1\cr
\vdots\cr p_n\end{pmatrix}\rangle = (-1)^n \langle u,\hbox{\rm
div}^{\{n\}} \begin{pmatrix}p_0\cr p_1\cr \vdots\cr
p_n\end{pmatrix}\rangle =
(-1)^n \sum_{i=0}^n \langle u, \mathcal{D}_i^n \, p_i\rangle, \\
&~& \\
&~&  \langle \hbox{\rm div}^{\{n\}} \begin{pmatrix}u_0\cr u_1\cr
\vdots\cr u_n\end{pmatrix}, p \rangle = (-1)^n \langle
\begin{pmatrix}u_0\cr u_1\cr \vdots\cr u_n\end{pmatrix}, \nabla^{\{n\}}  p
\rangle = (-1)^n \sum_{i=0}^n \langle u_i, \mathcal{D}_i^n\,
p\rangle.
\end{eqnarray*}

\section{Classical orthogonal polynomials in two variables}

In \cite{FPP1,FPP2,FPP3,FPP4}, the authors extended the concept of
{\it classical} bivariate orthogonal polynomials. In fact, we define
{\it classical} orthogonal polynomial starting from a matrix partial
differential equation with matrix coefficients, a direct
generalization of the partial differential equation studied by H. L.
Krall and I. M. Sheffer in \cite{KS}, and P. K. Suetin in \cite{Su}.

Let $a(x,y)$, $b(x,y)$, $c(x,y)$ be polynomials in two variables of total degree less
than or equal to 2, and let $d(x,y)$, $e(x,y)$ be polynomials in two variables of total degree less
than or equal to 1. Define the partial differential operator $L$ acting over
${\mathcal P}$ by means of
$$
L[p] \equiv a
\partial_{x x}p + 2 b \partial_{x y} p + c\partial_{y y}p + d \partial_x p + e \partial_y
p.
$$
The operator $L[\cdot]$ preserves the degree of the polynomials, that is,
$L[\mathcal{P}_n]\subset {\mathcal P}_n$.

Observe that $L$ is the left hand side of the partial differential
equation (\ref{KSorig}) studied in \cite{KS} and \cite{Su},
without restrictions on the polynomial coefficients.

Define
\begin{equation}\label{phipsi}
    \Phi = \begin{pmatrix}a& b\cr b&c\end{pmatrix} \in {\mathcal M}_{2}({\mathcal
                P_2}),\qquad
    \Psi =\begin{pmatrix}d\cr e\end{pmatrix} \in {\mathcal M}_{2\times 1}({\mathcal
            P_1}).
\end{equation}
Then, we can write
$$
L[p]\equiv\hbox{\rm div~}(\Phi\nabla p)+\tilde{\Psi}^t\nabla p,
$$
where $\tilde{\Psi}=\Psi - (\hbox{\rm div~} \Phi)^t$.

\begin{defn}
Let $u$ be a quasi definite moment functional, and let
$\{\mathbb{P}_n\}_{n\ge0}$ be the monic WOPS with respect to $u$.
Then, we say that $u$ is \emph{classical} (or
$\{\mathbb{P}_n\}_{n\ge0}$ is a \emph{classical} WOPS), if there
exist non--singular matrices $\Lambda_n\in{\mathcal
M}_{n+1}(\mathbb{R})$, such that, for $n\ge 1$,
\begin{equation}\label{de2}
L[\mathbb{P}^t_n]\equiv\hbox{\rm div~}(\Phi\nabla
\mathbb{P}^t_n)+\tilde{\Psi}^t\nabla\mathbb{P}^t_n =
\mathbb{P}^t_n \Lambda_n,
\end{equation}
and
\begin{equation}\label{phireg}
\det \, \langle u, \Phi\rangle \neq 0.
\end{equation}
\end{defn}

Differential equation (\ref{de2}) can be write also for $n=0$,
taking $\Lambda_0=0$.

Observe that, for $n=1$, then $\nabla \mathbb{P}_1^t =I_2$, and
equation (\ref{de2}) can be written as
$$\hbox{\rm div~}\Phi + \tilde{\Psi}^t = \mathbb{P}^t_1
\Lambda_1 \quad \Longrightarrow \quad \Psi^t = \mathbb{P}^t_1
\Lambda_1.$$ The non--singular character of $\Lambda_1$ implies
that $d$ and $e$ are independent polynomials of exact degree 1.

In \cite{KS} and \cite{Su}, every orthogonal polynomial of total
degree $n$ satisfies the same partial differential equation, and
therefore $\Lambda_n = \lambda_n \, I_{n+1}$, is a scalar matrix.

Let $L^*$ be the {\it formal Lagrange adjoint} of $L$, defined by
means of
\begin{equation}\label{lstar}
L^*[u]\equiv\hbox{\rm div~}(\Phi\nabla u) - \hbox{\rm
div~}(\tilde{\Psi}u),
\end{equation}
then it satisfies \quad $ \langle L^*[u], p\rangle = \langle u, L[p]\rangle,
\quad \forall p\in {\mathcal P}.
$

Observe that, if $u$ is classical, then $L^*[u] =0$. In fact, for any $n\ge 0$,
$$\langle L^*[u], \mathbb{P}_n^t \rangle = \langle u,
L[\mathbb{P}_n^t]\rangle = \langle u, \mathbb{P}_n^t\rangle\,
\Lambda_n =0,$$ since $\Lambda_0 = 0$, and $\langle u,
\mathbb{P}_n^t\rangle =0$, for $n\ge 1$, using the orthogonality of
the polynomials.

In \cite{FPP2}, the authors obtained several characterizations for
bivariate classical orthogonal polynomials. In particular, it is
shown that, under some mild regularity conditions, a quasi definite
moment functional $u$ is classical if and only if $u$ satisfies the
\emph{matrix Pearson--equation}
\begin{equation}\label{Pear}
\hbox{\rm div~} (\Phi \, u) = \Psi^t \, u,
\end{equation}
where $\Phi$ and $\Psi$ are the same polynomial matrices defined in (\ref{phipsi}).
Equivalently, $u$ is classical if and only
if
\begin{equation}\label{Pear2}
 \Phi \, \nabla ~u = \tilde \Psi ~u.
\end{equation}

\section{The symmetry factor}

As L. L. Littlejohn did in \cite{Li}, in order to obtain a
Rodrigues formula, we consider symmetry factors for the
differential operator $L[\cdot]$.

\begin{defn}[\cite{Li}] We say that $L[\cdot]$ is symmetric if $L[\cdot] = L^*[\cdot]$. $L[\cdot]$
is symmetrizable if there exists a
 nontrivial function $\omega(x,y)$ such that it is
$C^2$ in some open set, and $\omega \, L[\cdot]$ is symmetric. In
this case, $\omega$ is called a symmetry factor for $L[\cdot]$.
\end{defn}

As it is well known, if $\omega$ is a symmetry factor for
$L[\cdot]$, then the operator $\omega L$ can be written in a
compact form,
\begin{equation}\label{s-edgr}
\omega L[\centerdot] = {\hbox{\rm div~}}[\omega\,
\Phi\,\nabla\,\centerdot].
\end{equation}

From this equation, we can show that a function $\omega$ is a
symmetry factor for $L[\cdot]$ if and only if $\omega$ satisfies
the matrix Pearson--type equation (\ref{Pear2}).

\begin{proposition}
Let $\omega$ be a nontrivial function.  Then, $\omega$ is a
symmetry factor for equation (\ref{de2}), if and only if $\omega$
satisfies the matrix Pearson--type equation
\begin{equation}\label{pearm}
\Phi\,\nabla\,\omega = \tilde{\Psi}\,\omega.
\end{equation}
\end{proposition}

\begin{proof}
From equation (\ref{de2}), we get
\begin{eqnarray*}
 \omega \, L[\mathbb{P}_n^t]  &=& \omega [{\hbox{\rm
div}}(\Phi\,\nabla \mathbb{P}^t_n) +
\tilde{\Psi}^t\,\nabla\,\mathbb{P}^t_n ] = \\
&=& {\hbox{\rm div}}[\omega \, \Phi\,\nabla\,\mathbb{P}^t_n] -
(\nabla\,\omega)^t \Phi\,\nabla\,\mathbb{P}^t_n + \omega
\tilde{\Psi}^t\,\nabla\,\mathbb{P}^t_n =\\
&=&{\hbox{\rm div}}[\omega\,\Phi\,\nabla\,\mathbb{P}^t_n] -
 (\Phi\,\nabla\,\omega - \tilde{\Psi}\,\omega)^t\,\nabla\,\mathbb{P}^t_n,
    \end{eqnarray*}
   and the result follows using (\ref{s-edgr}).
\end{proof}

Observe that (\ref{pearm}) is equivalent to
\begin{equation}\label{tpms}
\hbox{\rm div~} (\Phi ~\omega) = \Psi^t ~\omega.
\end{equation}

The explicit expression for equations (\ref{pearm}) and
(\ref{tpms}) are
\begin{equation}\label{tpms2}
\left.\begin{array}{lll}
(a \, \omega)_x + (b \,\omega)_y &=& d\, \omega,\\
(b \, \omega)_x + (c \,\omega)_y &=& e\, \omega,
\end{array}\right\}\quad  \Longleftrightarrow\quad \left\{\begin{array}{lll}
a \, \omega_x + b \,\omega_y &=& (d-a_x - b_y) \, \omega,\\
b \, \omega_x + c \,\omega_y &=& (e - b_x - c_y)\, \omega.
\end{array}\right.
\end{equation}

 Obviously, the nontrivial solutions of the above system of partial differential equations give
us the symmetry factors for $L[\cdot]$.

\begin{proposition}[\cite{KKL2}] Suposse that $a\,c - b^2\not \equiv 0$. Then, the differential operator $L[\cdot]$
is symmetrizable if and only if
$$\frac{\partial}{\partial_y}\left[\frac{c(d-a_x-b_y) - b(e-b_x-c_y)}{a\,c - b^2}\right] =
\frac{\partial}{\partial_x}\left[\frac{a(e-b_x-c_y)-b(d-a_x-b_y)
}{a\,c - b^2}\right].
$$
\end{proposition}

In \cite{KKL2}, the authors shows that in the Krall and Sheffer
case, the existence of a symmetry factor is a necessary condition
for the existence of an OPS solution of the PDE (\ref{de2}). A quite
similar proof shows that the result is also true in the general
case.

\begin{proposition} If the
differential equation (\ref{de2}) has an OPS as solutions, and
$a\,c - b^2\not \equiv 0$, then $L[\cdot]$ must be symmetrizable.
\end{proposition}

\section{The second kind Kronecker power of a matrix}

The {\bf second kind Kronecker power} is defined in \cite{Be} p.
236, for a square matrix $A=(a_{i,j})$ of dimension 2. Let $z =
A\, t$ be the linear transformation defined by $A$, that is,
\begin{equation}\label{SKP}\left(%
\begin{array}{c}
  z_1 \\
  z_2 \\
\end{array}%
\right) =\left(%
\begin{array}{cc}
  a_{0,0} & a_{0,1} \\
  a_{1,0} & a_{1,1} \\
\end{array}%
\right)\left(%
\begin{array}{c}
  t_1 \\
  t_2 \\
\end{array}%
\right) \Rightarrow \left\{\begin{array}{ccccc}
                        z_1 & = & a_{0,0}\,t_1 & + &
                          a_{0,1} \,t_2, \\
                          z_2 & = & a_{1,0} \,t_1 & + &
                          a_{1,1}\,t_2.\\
                        \end{array}\right.
\end{equation}
 Each of the $(n+1)$ homogeneous products
$z_1^{n-i} z_2^{i}$, $i=0,1,\ldots,n$, is transformed under
(\ref{SKP}) into a linear combination of the $(n+1)$ homogeneous
products $t_1^{n-i} t_2^{i}$, $i=0,1,\ldots,n$. Then, the square
matrix of order $(n+1)$ specifying these linear transformation is
known as the {\bf second kind n--th Kronecker power} of $A$, which
we denote by $A^{\{n\}}$.

 Throughout, we will adopt the convenion $A^{\{0\}}\equiv
1$. Observe that $A^{\{1\}} = A$. Moreover, we can check that
$$
A^{\{2\}} = \left(
 \begin{array}{ccc} a_{0,0}^2 & 2\, a_{0,0}\, a_{0,1}&
a_{0,1}^2\\ a_{0,0}\, a_{1,0}& (a_{0,0}\, a_{1,1} + a_{0,1}\,
a_{1,0}) &
        a_{0,1}\, a_{1,1}\\
        a_{1,0}^2 & 2\, a_{1,0}\, a_{1,1} & a_{1,1}^2
        \end{array}\right).
        $$
It is possible to give the explicit expression for the entries of
$A^{\{n\}}$, $n\ge 0$. In fact, if we denote $A^{\{n\}} =
(a^{\{n\}}_{i,j})_{i,j=0}^n$, then, a direct calculation shows that
\begin{equation}\label{skpex}
a_{i,j}^{\{n\}} = \sum_{k=0}^j \,\binom{n-i} {k}\,\binom{i}{j-k}\,
a_{0,0}^{n-i-k} \, a_{0,1}^{k} \, a_{1,0}^{i-j+k} \,
a_{1,1}^{j-k}, \quad 0\le i,j\le n,
\end{equation}
where, as usual, $\binom{m}{l} = 0$, if $m<l$.

From the above explicit expression we can obtain recurrence
formulas for the second kind Kronecker power of a matrix $A$.

\begin{lemma}\label{SKPrecu}
There exist two recurrence formulas for $A^{\{n\}}$, $n\ge 1$. In
fact, for $0\le j\le n,$ we get

\noindent \emph{Recurrence I:}
$$\left\{\begin{array}{lll}
a_{i,j}^{\{n\}} &=& a_{0,0} \, a_{i,j}^{\{n-1\}} + a_{0,1}\,a_{i,j-1}^{\{n-1\}},\quad 0\le i\le n-1,\\
a_{n,j}^{\{n\}} &=& a_{1,0} \, a_{n-1,j}^{\{n-1\}} +
a_{1,1}\,a_{n-1,j-1}^{\{n-1\}}.
\end{array}\right.$$
\noindent \emph{Recurrence II:}
$$\left\{\begin{array}{lll}
a_{0,j}^{\{n\}} &=& a_{0,0} \, a_{0,j}^{\{n-1\}} + a_{0,1}\,a_{0,j-1}^{\{n-1\}},\\
a_{i,j}^{\{n\}} &=& a_{1,0} \, a_{i-1,j}^{\{n-1\}} +
a_{1,1}\,a_{i-1,j-1}^{\{n-1\}}, \quad 1\le i\le n.
\end{array}\right.$$

In a matrix form, we can express the above formulas as follows

\noindent \emph{Recurrence I:}
$$
A^{\{n\}} = \left(%
\begin{array}{ccc}
  a_{0,0} &  & \bigcirc   \\
         &  \ddots &  \\
  \bigcirc   & & a_{0,0} \\
  \hline
    0 &   \cdots & a_{1,0}  \\
\end{array}%
\right)A^{\{n-1\}} L_{n-1}^0 + \left(%
\begin{array}{ccc}
  a_{0,1} & & \bigcirc   \\
         & \ddots &  \\
  \bigcirc  & & a_{0,1} \\
  \hline
    0 &  \cdots & a_{1,1}  \\
\end{array}
\right)A^{\{n-1\}} L_{n-1}^1.
$$
\noindent \emph{Recurrence II:}
$$
A^{\{n\}} = \left(%
\begin{array}{ccc}
  a_{0,0} & \cdots & 0  \\
  \hline
  a_{1,0}  & & \bigcirc   \\
          & \ddots &  \\
  \bigcirc  & & a_{1,0} \\
\end{array}%
\right)A^{\{n-1\}} L_{n-1}^0 + \left(%
\begin{array}{cccc}
  a_{0,1} & \cdots & 0  \\
  \hline
  a_{1,1} & & \bigcirc   \\
          & \ddots &  \\
  \bigcirc  & & a_{1,1} \\
\end{array}%
\right)A^{\{n-1\}} L_{n-1}^1,
$$
where $L_{n-1}^k$, $k=0,1$, are $n\times(n+1)$ matrices defined as
\begin{equation}\label{li}
L_{n-1}^0 = \left(%
\begin{array}{cccc}
  1 & ~         & \bigcirc   & 0 \\
  ~ & \ddots    & ~          & \vdots \\
  \bigcirc & ~  & 1          & 0 \\
\end{array}%
\right), \quad \hbox{\rm and} \quad L_{n-1}^1 = \left(%
\begin{array}{cccc}
  0 &  1 & ~& \bigcirc \\
  \vdots & ~ & \ddots & ~ \\
  0 & \bigcirc & ~ & 1 \\
\end{array}%
\right).\end{equation}
\end{lemma}

\section{The matrix Rodrigues--type formula}\label{mainsec}

In this Section, we shall assume that the differential operator
$L[\cdot]$ is symmetrizable, and we denote by $\omega$ a
nontrivial symmetry factor of $L[\cdot]$. Our main result is
stated in the following theorem.

\begin{theorem}\label{main} Let $u$ be a classical moment functional, and let (\ref{de2}) be the
matrix partial differential equation associated with $u$. Let
$\omega$ be a nontrivial symmetry factor of $L[\cdot]$, and let us
assume that there exist polynomial matrices $\Psi_0, \Psi_1\in
\mathcal{M}_{2}(\mathcal{P}_1)$, such that
\begin{eqnarray}
&~& (a \, \Phi)_x + (b \, \Phi)_y  =  \Phi\, \Psi_0,\nonumber\\
&~& (b \, \Phi)_x + (c \, \Phi)_y = \Phi\, \Psi_1.\label{condR}
\end{eqnarray}

Then, for $n\ge0$, the expression
\begin{equation}\label{rod}
\mathbb{Q}_n^t = \frac{1}{~\omega~}\, {\hbox{\rm div}}^{\{n\}}
(\Phi^{\{n\}}\,\omega),\qquad n\ge 0,
\end{equation}
provides a $1\times (n+1)$ polynomial vector of degree
 $n$, such that
$$\langle u,\mathbb{Q}_n \, p\rangle = 0, \quad \forall p\in \mathcal{P}_{n-1}.$$
Moreover, if $\{\mathbb{Q}_n\}_n$ is a PS, then it is a WOPS
associated with $u$.
\end{theorem}

As we will show later, we can establish sufficient conditions in
order to obtain WOPS from the matrix Rodrigues formula.

\begin{cor}\label{cordiag} In the hypothesis of Theorem \ref{main}, if $\Phi$ is a diagonal matrix, then
$\{\mathbb{Q}_n\}_n$ is a WOPS associated with $u$.
\end{cor}

\begin{cor}\label{cor2} In the hypothesis of Theorem \ref{main}, if the $(n+1)$ square matrix
$$\langle u, \Phi^{\{n\}}\rangle,$$
is non--singular, $n\ge0$, then $\{\mathbb{Q}_n\}_n$ is a WOPS
associated with $u$.
\end{cor}

Next, we will study some particular cases. Moreover, we will deduce
that the Rodrigues--type formula given by P. K. Suetin (\cite{Su})
and revisited in Y. J. Kim et al. (\cite{KKL2}) is a particular case
of the matrix Rodrigues--type formula (\ref{rod}).

\bigskip

\noindent{\bf Example 1: The diagonal case}

Let us assume that $\Phi$ is diagonal, i. e., $b\equiv 0$. In this
case, the matrix Pearson--type equation (\ref{tpms2}) reduces to
\begin{equation}\label{diag}
\left(a\,\omega \right)_x = d \, \omega,\qquad \left(c\,\omega
\right)_y = e \, \omega,
\end{equation}
and the $2\times 2$ matrices $\Psi_k = (\psi_{i,j}^k)_{i,j=0}^1$,
$k=0,1$, in (\ref{condR}) are polynomials of degree less than or
equal to 1 satisfying
$$(a \, \Phi)_x  =  \Phi\, \Psi_0,\qquad  (c \, \Phi)_y = \Phi\,
\Psi_1.
$$
This condition holds if we take $\psi_{0,1}^k=\psi_{1,0}^k = 0$,
$\psi_{0,0}^0=2\,a_x$, $\psi_{1,1}^1 =2\,c_y$, and
$\psi_{1,1}^{0}$, $\psi_{0,0}^{1}$,  such that
$$a_x\,c + a\,c_x = c\, \psi_{1,1}^{0},\qquad a_y\,c + a\,c_y = a\, \psi_{0,0}^{1}.
$$

In this way, if, for instance, the expressions
\begin{equation}\label{cond}
\frac{a\, c_x}{c}, \qquad {\hbox{\rm and}}\qquad
\frac{a_y\,c}{a},\end{equation} are polynomials of total degree
less than or equal to 1, then (\ref{condR}) holds. This is
possible, for example, when $a$ and $c$ are equal up to a
multiplicative constant factor, or $a_y = c_x = 0$, that is, the
polynomial $a$ only depends on $x$, and $c$ only depends on $y$.
This last situation corresponds to the case where the moment
functional is a tensor product of univariate moment functionals.

Moreover, the second kind n--th order Kronecker power of $\Phi$,
is again a diagonal matrix whose elements are given by
$$\phi_{i,i}^{\{n\}} = a^{n-i} \, c^{i}, \qquad i=0,1,\ldots, n,$$
and the Rodrigues formula gives
$$\mathbb{Q}_n^t = (Q_{n,0}, Q_{n-1,1}, \ldots, Q_{0,n})^t,
$$
where
$$Q_{n-i,i}(x,y) = \frac{1}{~\omega~}\binom{n}{i}\, \partial_x^{n-i}\,\partial_y^{i}\,
(a^{n-i}\, c^{i}\,\omega), \quad 0\le i \le n.
$$
In addition, from Corollary \ref{cordiag}, we deduce that
$\{\mathbb{Q}_n\}_n$ is a WOPS associated to $u$.

\bigskip

\noindent{\bf Example 2: Tensor product of classical orthogonal
polynomials in one variable}

Tensor product of two families of classical orthogonal polynomials in one
variable (Hermite,
Laguerre, Jacobi or Bessel), $\{R_h\}_{h\ge0}$ and $\{S_k\}_{k\ge0}$ is defined by means of
$$
P_{h,k}(x,y)=R_{h}(x) S_{k}(y), \quad h, k \ge 0.
$$
These families of two--dimensional polynomials are classical
according to our definition (see \cite{FPP2}). The symmetry factor
for the partial differential equation (\ref{de2}) is given by
$$
\omega(x,y)=\omega_1(x) \, \omega_2(y),
$$
where $\omega_1(x)$ and $\omega_2(y)$ are symmetry factors of the
differential equations for the polynomials $\{R_h(x)\}_{h\ge0}$ and
$\{S_k(y)\}_{k\ge0}$, respectively. This is a particular case of Example 1, since
$$\Phi = \begin{pmatrix}a & 0 \cr
                0 & c\end{pmatrix},\qquad
    \Psi =\begin{pmatrix}d \cr e\end{pmatrix},
$$
and
$$(a \, \omega_1)_x = d \, \omega_1,\qquad (c \,
\omega_2)_y = e \, \omega_2,$$ are the respective Pearson equations
for the univariate polynomials. Observe that $a_y = d_y =0$, that
is, $a$ and $d$ only depend on $x$; and $c_x = e_x =0$, that is, $c$
and $e$ only depend on $y$. Then, (\ref{condR}) holds and the
Rodrigues formula gives
$$\mathbb{Q}_n^t = (Q_{n,0}, Q_{n-1,1}, \ldots, Q_{0,n})^t,
$$
where
$$Q_{n-i,i}(x,y) = \binom{n}{i} \frac{1}{~\omega_1~}
\partial_x^{n-i} \,(a^{n-i}\,\omega_1) \,\frac{1}{~\omega_2~} \, \partial_y^{i}\,
(c^{i}\,\omega_2) , \quad 0\le i \le n.
$$
Observe that this expression provides $Q_{n-i,i}(x,y)$ as the
product of the one dimensional Rodrigues formulas for
$\{R_h\}_{h\ge0}$ and $\{S_k\}_{k\ge0}$ (see, for example,
\cite{Chi}), up to a constant factor.

General tensor product of univariate classical orthogonal
polynomials is not included in the Krall and Sheffer (\cite{KS}),
and Suetin (\cite{Su}) classifications for the bivariate classical
orthogonal polynomials, except for Hermite and Laguerre families.

\bigskip

\noindent{\bf Example 3: The Suetin--Rodrigues formula}

Under several conditions, P. K. Suetin (\cite{Su}) and Y. J. Kim
et al (\cite{KKL2}) proved a Rodrigues formula for some classical
orthogonal polynomials in two variables in the Krall and Sheffer
sense. We will show that this Rodrigues formula is again a
particular case of the matrix Rodrigues formula given in
(\ref{rod}).

By elementary algebraic manipulations, equation (\ref{tpms2}) can
be reduces to the form
$$\left(%
\begin{array}{cc}
  \alpha & 0 \\
  0 &\alpha \\
\end{array}%
\right)\, \, \left(%
\begin{array}{c}
  \omega_x \\
  \omega_y \\
\end{array}%
\right) \, = \, \left(%
\begin{array}{c}
  \beta \\
  \gamma\\
\end{array}%
\right)\, \omega
\Longrightarrow \left\{\begin{array}{c}
                       \alpha \, \omega_x = \beta\, \omega,\\
\alpha \, \omega_y = \gamma\, \omega,\end{array}\right.
$$
where $\alpha = \det \Phi$, $\beta = c(d-a_x-b_y) - b(e-b_x-c_y)$, and $\gamma = -b(d-a_x-b_y) + a(e-b_x-c_y).$

In the Suetin case, the special elections for the polynomials $a$,
$b$, $c$, $d$ and $e$, makes $\deg(\alpha)\le 3$, $\deg(\beta)\le
2$, and $\deg(\gamma)\le 2$. Then, if we assume
$$ a = a_1\,a_2, \quad b = a_1\, b_1\, c_1, \quad c = c_1\,c_2,
$$ where $a_1\neq 0$, $c_1\neq0$. Then
$$\alpha = a_1\, c_1 \, \alpha_0,\qquad \beta = c_1 \beta_0,\qquad \gamma = a_1\, \gamma_0,$$
and $\omega$ satisfies the following system of partial
differential equations
$$
p \, \omega_x = \beta_0\, \omega, \qquad q \, \omega_y = \gamma_0\, \omega,
$$
where $p=a_1\,\alpha_0$, and $q=c_1\,\alpha_0$.

In \cite{KKL2} and \cite{Su}, the authors prove that if $(a_1)_y =
(c_1)_x =0$, $\deg(p),\,\deg(q)\le 2$, and $\deg(\beta_0),\,
\deg(\gamma_0)\le 1$, then the sequence $\{\mathbb{P}_n\}_n$, with
$$\mathbb{P}_n=(P_{n,0},P_{n-1,1},\ldots,P_{0,n})^t,$$ defined by
means of
$$P_{n-i,i}(x,y) = \frac{1}{~\omega~} \,\partial^{n-i}_x\,\partial_y^{i}(p^{n-i}\, q^i\,\omega),$$
provides a WOPS satisfying the partial differential equation
(\ref{KSorig}).

This is again a particular case of Example 1, since Suetin hypothesis $(a_1)_y = (c_1)_x =0$ implies
(\ref{condR}). Matrix Rodrigues formula (\ref{rod}) provides, for $n\ge0$, a WOPS
$$\mathbb{Q}_n^t = (Q_{n,0}, Q_{n-1,1}, \ldots, Q_{0,n})^t,
$$
where
$$Q_{n-i,i}(x,y) = \frac{1}{~\omega~}\binom{n}{i}\, \partial_x^{n-i}\,\partial_y^{i}\,
(p^{n-i}\, q^{i}\,\omega), \quad 0\le i \le n.
$$
Observe that the above expression coincides with Suetin--Rodrigues
formula up to the binomial coefficients.

\bigskip

\noindent{\bf Example 4: Orthogonal polynomials on the unit ball}

Classical orthogonal polynomials on the unit ball, $B_2=\{(x,y):
x^2 + y^2\leq 1\}$, are associated with the weight function (symmetry factor)
$$
\omega(x,y)=(1-x^2-y^2)^{\mu-1/2}, \qquad \mu
>-1/2.
$$
In this case, the matrices $\Phi$ and $\Psi$ are given by
$$\Phi = \begin{pmatrix}
                x^2-1 & x\, y\cr
                x\, y & y^2-1\end{pmatrix},\qquad
    \Psi =\begin{pmatrix}(2\mu +2)x \cr (2\mu +2)y\end{pmatrix},
$$
and the matrix Pearson--type equation for the symmetry factor $\omega$ is given by
$$\begin{pmatrix} x^2-1 & x\, y\cr
                x\, y & y^2-1\end{pmatrix}\, \begin{pmatrix}\omega_x\cr \omega_y\end{pmatrix} \, =\,
    \begin{pmatrix}(2\mu -1)x \cr (2\mu -1)y\end{pmatrix} \, \omega.
    $$
We can check that condition (\ref{condR}) is not true, and matrix
Rodrigues formula (\ref{rod}) does not provides a WOPS. However,
we can transform the above matrix Pearson--type equation and
obtain
\begin{equation}\label{circle}
\begin{pmatrix}1 - x^2- y^2 & 0\cr
                0 & 1 - x^2 - y^2\end{pmatrix}\, \begin{pmatrix}\omega_x\cr \omega_y\end{pmatrix} \, =\,
    \begin{pmatrix}-(2\mu -1)x \cr -(2\mu -1)y\end{pmatrix} \, \omega.
\end{equation}
Equation (\ref{circle}) is a Pearson--type equation for the
symmetry factor $\omega$, it is diagonal, and it satisfies
condition (\ref{cond}). Then, matrix Rodrigues formula provides a
WOPS defined as follows
$$\mathbb{Q}_n^t = (Q_{n,0}, Q_{n-1,1}, \ldots, Q_{0,n})^t,
$$
where \begin{equation}\label{rodbola}Q_{n-i,i}(x,y) =
\frac{1}{~\omega~}\binom{n}{i}\,
\partial_x^{n-i}\,\partial_y^{i}\, ((1-x^2-y^2)^{n}\,\omega),
\quad 0\le i \le n,
\end{equation}
which coincides with the Rodrigues formula for the orthogonal
polynomials on the unit ball as described in \cite{DX,KKL2,Su}.

\bigskip

\noindent{\bf Example 5: Orthogonal polynomials on the simplex
(Appell polynomials)}

The weight function associated with the classical orthogonal
polynomials on the simplex, $T=\{(x,y):\,x,y\geq 0, 1-x-y\geq
0\}$, is defined by
$$
\omega(x,y)=x^{\alpha}y^{\beta}(1-x-y)^{\gamma}, \qquad \alpha,
\beta, \gamma > -1.
$$
Here, the matrices $\Phi$ and $\Psi$ are given by
$$\Phi = \begin{pmatrix}
                x(x-1) & x y \cr
                x y & y(y-1) \end{pmatrix}, \qquad
    \Psi =\begin{pmatrix}
           (\alpha+\beta+\gamma+3)x-(\alpha+1) \cr
           (\alpha+\beta+\gamma+3)y-(\beta+1)\end{pmatrix}.$$
Observe that condition (\ref{condR}) holds, and therefore, the
matrix Rodrigues formula for triangle polynomials provides a WOPS.
Since the obtained polynomials are monic, they coincide with the
so--called {\it second kind Appell polynomials} (see \cite{Su}).

Y. J. Kim et al (see \cite{KKL2}) proved that $\omega$ satisfies
also a diagonal matrix Pearson--type equation
$$
\begin{pmatrix}x\,(1-x-y) & 0\cr
                0 & y\,(1-x-y)\end{pmatrix}\, \begin{pmatrix}\omega_x\cr \omega_y\end{pmatrix} \, =\,
    \begin{pmatrix}\alpha(1-x-y)-\gamma\,x \cr \beta (1-x-y)-\gamma\,y\end{pmatrix} \, \omega.
$$
In this way, condition (\ref{condR}) holds, and we are again in
the situation described in Example 1. Then, we get
$$\mathbb{Q}_n^t = (Q_{n,0}, Q_{n-1,1}, \ldots, Q_{0,n})^t,
$$
defined from the matrix Rodrigues formula
\begin{equation}\label{rodtrian}
Q_{n-i,i}(x,y) = \frac{1}{~\omega~}\binom{n}{i}\,
\partial_x^{n-i}\,\partial_y^{i}\, (x^{n-i}\, y^{i}
(1-x-y)^n\,\omega), \quad 0\le i \le n, \end{equation} is a WOPS
relative to $\omega$.

\bigskip

Rodrigues formula for the circle and the triangle polynomials,
(\ref{rodbola}) and (\ref{rodtrian}) respectively, coincide with the
classical expressions for these polynomials which appear in the
literature (see, for instance, \cite{AK,DX,KKL2,Koor,Su}).

\bigskip

\noindent{\bf Example 6: The most intriguing case (sic, L. L. Littlejohn, \cite{Li})}

Krall and Sheffer (\cite{KS}) showed
that the differential equation
$$L[p] \equiv 3 y p_{xx} + 2 p_{x y} -x p_x-y p_y = -n p,$$
has an OPS as solutions. In this case, the matrices $\Phi$ and
$\Psi$ are given by
$$\Phi = \begin{pmatrix}3y & 1 \cr
                1 & 0 \end{pmatrix}, \qquad
    \Psi =\begin{pmatrix}-x \cr -y\end{pmatrix}.$$
A symmetry factor for the partial differential equation is
$\omega(x,y) = \exp(y^3-xy)$ (\cite{Li}). Observe that condition
(\ref{condR}) is satisfied, and then we can obtain a matrix
Rodrigues formula for these polynomials. In particular, we get
\begin{eqnarray*}
\mathbb{Q}_0^t &=& 1, \\
\mathbb{Q}_1^t &=& \left(-x, -y\right), \\
\mathbb{Q}_2^t &=& \left(x^2 - 6 y, 2 x y-2 , y^2\right), \\
\mathbb{Q}_3^t &=& \left(-x^3+18xy-12,-3x^2y +18y^2 +
6x,-3xy^2+6y,-y^3 \right).
\end{eqnarray*}

We must remark that Rodrigues--type formula for these polynomials
can not be obtained using the Suetin tools (see \cite{KKL2}).
Moreover, the symmetry factor $\omega$ is not a weight function.
These polynomials has attracted considerable attention
(\cite{KKL1,KKL2,KS,Li,Su}) since this is the simplest case of
classical orthogonal polynomials in two variables orthogonal with
respect to a non positive definite moment functional.

\section{Proof of the main result}

In this section, we will prove the matrix Rodrigues formula,
introduced in Theorem \ref{main}, in several steps that we will
organize in a series of lemmas.

\begin{lemma}\label{lemacondR}
Let $A=(a_{i,j})_{i,j=0}^1$ be a $2\times 2$ polynomial matrix.
Assume that, there exist polynomial matrices $\Psi_k =
(\psi^k_{i,j})_{i,j=0}^1$,
 such that
\begin{equation}
(a_{k,0} \, A)_x + (a_{k,1} \, A)_y = A \,\Psi_k,\qquad k=0,1.
\end{equation}

Then, for $n\ge 1$, we have
\begin{equation}\label{condRAn}
(a_{k,0} \, A^{\{n\}})_x + (a_{k,1} \, A^{\{n\}})_y = A^{\{n\}}\,
\Psi_k^n,\qquad k= 0,1, \end{equation} where
$\Psi_k^n=(\psi_{i,j}^{n,k})_{i,j=0}^n\in
\mathcal{M}_{n+1}(\mathcal{P})$, $k=0,1$, are three--diagonal
matrices with
$$
\begin{array}{lclll}
\psi_{j-1,j}^{n,k} &=& (n+1-j)\, \psi_{0,1}^k, & & 1\le j
\le n,\\
\psi_{j,j}^{n,k} &=& (n-j)\, \psi_{0,0}^k + j \, \psi_{1,1}^k -
(n-1)\left[(a_{k,0})_x + (a_{k,1})_y\right], &
& 0\le j \le n,\\
\psi_{j+1,j}^{n,k} &=& (j+1)\, \psi_{1,0}^k, & & 0\le j \le n-1.
\end{array}
$$
Moreover, $\deg\, \Psi_k^n \le \deg\, A - 1$, \, $k=0, 1$, \,
$n\ge 1$.
\end{lemma}

\begin{proof}
The Lemma follows by induction on $n$. For $n=1$, the result holds
using $A^{\{1\}}=A$, and $\Psi_k^1 = \Psi_k$, $k=0,1$.

Let $0\le i\le n-1$, and $0\le j\le n$, and assume that the result
holds for $n-1$. First, we compute the left hand side of
(\ref{condRAn}) using the induction hypothesis, and Recurrence I
of Lemma \ref{SKPrecu}.

\begin{eqnarray*}
\lefteqn{\left(a_{0,0} \, a_{i,j}^{\{n\}}\right)_x +
\left(a_{0,1}\,
a_{i,j}^{\{n\}}\right)_y=}\\
 &=& \left(a_{0,0}\left[a_{0,0} \, a_{i,j}^{\{n-1\}} + a_{0,1} a_{i,j-1}^{\{n-1\}}\right]\right)_x +
 \left(a_{0,1}\,
\left[a_{0,0} \, a_{i,j}^{\{n-1\}} + a_{0,1} a_{i,j-1}^{\{n-1\}}\right]\right)_y\\
&=& a_{0,0}\left[(a_{0,0}\, a_{i,j}^{\{n-1\}})_x +
(a_{0,1}\,a_{i,j}^{\{n-1\}})_y\right] + \left[(a_{0,0})_x\,
a_{0,0} + (a_{0,0})_y \, a_{0,1}\right]\,
a_{i,j}^{\{n-1\}} \\
&+& a_{0,1}\left[(a_{0,0}\, a_{i,j-1}^{\{n-1\}})_x +
(a_{0,1}\,a_{i,j-1}^{\{n-1\}})_y\right]+ \left[(a_{0,1})_x\,
a_{0,0} + (a_{0,1})_y \,
a_{0,1}\right]\, a_{i,j-1}^{\{n-1\}}\\
&=& a_{0,0} \sum_{l=j-1}^{j+1}a_{i,l}^{\{n-1\}}
\,\psi_{l,j}^{n-1,0} + \left\{a_{0,0}\,\left[\psi_{0,0}^{0} +
(a_{0,0})_x + (a_{0,1})_y\right] + a_{0,1} \,
\psi_{1,0}^0\right\}\, a_{i,j}^{\{n-1\}} \\
& + & a_{0,1}
\sum_{l=j-2}^{j}a_{i,l}^{\{n-1\}}\,\psi_{l,j-1}^{n-1,0} +
\left\{a_{0,0}\,\psi_{0,1}^0 + a_{0,1}\left[\psi_{1,1}^{0} +
(a_{0,0})_x + (a_{0,1})_y\right]
\right\}\, a_{i,j-1}^{\{n-1\}}\\
&=& a_{0,0}\, a_{i,j-1}^{\{n-1\}}\left[\psi_{0,1}^0 +
\psi_{j-1,j}^{n-1,0}\right] +
a_{0,1}\, a_{i,j-2}^{\{n-1\}} \, \psi_{j-2,j-1}^{n-1,0} \\
&+& a_{0,0}\, a_{i,j}^{\{n-1\}}\left[\psi_{0,0}^0 +
\psi_{j,j}^{n-1,0} + (a_{0,0})_x + (a_{0,1})_y\right] \\
&+& a_{0,1}\, a_{i,j-1}^{\{n-1\}} \left[\psi_{1,1}^0 +
\psi_{j-1,j-1}^{n-1,0} + (a_{0,0})_x +
(a_{0,1})_y\right]\\
& +& a_{0,0}\, a_{i,j+1}^{\{n-1\}}\,\psi_{j+1,j}^{n-1,0} +
a_{0,1}\, a_{i,j}^{\{n-1\}}\left[\psi_{1,0}^0 +
\psi_{j,j-1}^{n-1,0}\right].
\end{eqnarray*}
Now, replace the recurrence relations for $\psi_{i}^{n-1,0}$, and
Recurrence I for the elements $a_{i,j}^{\{n\}}$ to obtain
$$\left(a_{0,0} \, a_{i,j}^{\{n\}}\right)_x + \left(a_{0,1}\,
a_{i,j}^{\{n\}}\right)_y= a_{i,j-1}^{\{n\}} \, \psi_{j-1,j}^{n,0}
+ a_{i,j}^{\{n\}}\,\psi_{j,j}^{n,0} + a_{i,j+1}^{\{n\}}\,
\psi_{j+1,j}^{n,0}.
$$
The rest of the cases follows in a similar way.
\end{proof}

\begin{lemma} Assume that the matrix polynomial $\Phi$ satisfies
(\ref{condR}). Let $A_n\in \mathcal{M}_{(n+1)\times
(m+1)}(\mathcal{P})$ be an arbitrary polynomial matrix. Then,
$${\hbox{\rm div}}^{\{n\}}
(\Phi^{\{n\}}\,A_n\,\omega) = \hbox{\rm div}^{\{n-1\}}
(\Phi^{\{n-1\}}\,A_{n-1}\,\omega),
$$
where $A_{n-1}\in \mathcal{M}_{n\times (m+1)}(\mathcal{P})$ is a
polynomial matrix satisfying
$$\deg A_{n-1} \le \deg A_n + 1.$$
\end{lemma}

\begin{proof} From the definition of ${\hbox{\rm div}}^{\{n\}}$ in (\ref{divn}),
we can easily deduce
$${\hbox{\rm div}}^{\{n\}} (\Phi^{\{n\}}\,A_n\,\omega) = {\hbox{\rm
div}}^{\{n-1\}} \left[\partial_x (L_{n-1}^0
\,\Phi^{\{n\}}\,A_n\,\omega) +
\partial_y (L_{n-1}^1\,\Phi^{\{n\}}\,A_n\,\omega)\right].$$
Now, if we apply recurrence formulas in Lemma \ref{SKPrecu}, we
get
\begin{eqnarray*}
\lefteqn{\partial_x (L_{n-1}^0 \,\Phi^{\{n\}}\,A_n\,\omega) +
\partial_y (L_{n-1}^1\,\Phi^{\{n\}}\,A_n\,\omega) =\quad \qquad \qquad}\\
&=& \left[(a\,\Phi^{\{n-1\}})_x + (b\, \Phi^{\{n-1\}})_y\right] L_{n-1}^0 A_n \,\omega \\
&~& +\left[(b\,\Phi^{\{n-1\}})_x + (c\, \Phi^{\{n-1\}})_y\right] L_{n-1}^1 A_n \,\omega\\
&~& + \Phi^{\{n-1\}} \left[(a\,L_{n-1}^0 + b \,L_{n-1}^1)(A_n)_x + (b\,L_{n-1}^0 + c \,L_{n-1}^1)(A_n)_y\right]\omega\\
&~& + \Phi^{\{n-1\}} \left[L_{n-1}^0\,A_n (a\,\omega_x +
b\,\omega_y) + L_{n-1}^1\,A_n (b\,\omega_x + c\,\omega_y)\right].
\end{eqnarray*}
Condition (\ref{condR}), Lemma \ref{lemacondR}, and equation
(\ref{pearm}) give
$${\hbox{\rm div}}^{\{n\}} (\Phi^{\{n\}}\,A_n\,\omega) = {\hbox{\rm div}}^{\{n-1\}} (\Phi^{\{n-1\}}\,A_{n-1}
\,\omega),$$ where
\begin{eqnarray*}
A_{n-1} &=& \left\{[\Psi_0^{n-1} + (d-a_x-b_y)I_n]L_{n-1}^0 + [\Psi_1^{n-1}+ (e-b_x-c_y)I_n] L_{n-1}^1\right\}A_n \\
&~& + (a\,L_{n-1}^0 + b\,L_{n-1}^1)(A_n)_x + (b\,L_{n-1}^0 +
c\,L_{n-1}^1)(A_n)_y.
\end{eqnarray*}
Observe that the degree condition can be deduced from the explicit
expression for $A_{n-1}$.
\end{proof}

Now, taking $A_n = I_{n+1}$, and applying induction on $n$, we get

\begin{lemma}\label{recu11} In the hypothesis of Lemma \ref{recu11}, for any $n\ge 0$,
$$\mathbb{Q}_n^t := \frac{1}{~\omega~}{\hbox{\rm div}}^{\{n\}}
(\Phi^{\{n\}}\,\omega),$$ is a $1\times(n+1)$ polynomial matrix of
degree less than or equal to $n$.
\end{lemma}

Using the same technique as Lemma \ref{recu11}, we can write
$\mathbb{Q}_n^t$ in terms of the moment functional $u$.

\begin{lemma} Let $u$ be a classical moment functional, and let (\ref{de2}) be the
matrix partial differential equation associated with $u$. Assume
that condition (\ref{condR}) holds. Then, for $n\ge 0$,
$\mathbb{Q}_n$ defined in Lemma \ref{recu11} satisfies

\noindent (i) $\mathbb{Q}_n^t \,u = {\hbox{\rm div}}^{\{n\}}
(\Phi^{\{n\}}\,u),\qquad n\ge 0,$

\bigskip

\noindent (ii) $\langle u, \mathbb{Q}_n \mathbb{X}_m^t\rangle =0,
\quad 0\le m\le n-1,$ where $\mathbb{X}_m^t =
(x^{m-i}\,y^i)_{i=0}^m$,

\bigskip

\noindent (iii) $\{\mathbb{Q}^t_n\}_n$ is a $1\times(n+1)$ vector of
polynomials of degree less than or equal to $n$. In fact, if we
denote
$$\mathbb{Q}_n^t = (Q_{n,0}, Q_{n-1,1}, \ldots,
Q_{0,n}),
$$
then $Q_{n-i,i}$ has exact degree $n$ or $Q_{n-i,i}\equiv 0$, for
$0\le i\le n$.

\end{lemma}

\begin{proof} (i) follows from a similar reasoning as used in Lemma \ref{recu11}.

In order to prove (ii), we compute
\begin{eqnarray*}
\langle u, \mathbb{Q}_n \mathbb{X}_m^t\rangle &=& \langle
\mathbb{Q}_n^t\, u, \mathbb{X}_m^t\rangle = \langle {\hbox{\rm
div}}^{\{n\}} (\Phi^{\{n\}}\,u), \mathbb{X}_m^t\rangle =\\
&=& (-1)^n \langle \,\Phi^{\{n\}}\,u,
\nabla^{\{n\}}\mathbb{X}_m^t\rangle =0, \quad 0\le m\le n-1.
\end{eqnarray*}
As a consequence, $\langle u, \mathbb{Q}_n\, p(x,y)\rangle =0$,
for any $p(x,y) \in \mathcal{P}_{n-1}$.

From the above property, we will prove (iii). Let $n\ge1$, and
suppose that there exists $0\le i\le n$, such that $\deg \,
Q_{n-i,i} <n$. Then,

$$\langle u, \mathbb{Q}_n \, Q_{n-i,i}\rangle =0,$$
and therefore $\langle u, Q_{n-i,i}^2\rangle =0$. Since $u$ is a
quasi definite moment functional, we deduce $Q_{n-i,i} \equiv 0$.

\end{proof}

\bigskip

To end the Section, we will prove Corollaries \ref{cordiag} and
\ref{cor2}.

\bigskip

\noindent{\bf Proof of Corollary \ref{cordiag}}

\begin{proof} If $\Phi$ is a diagonal matrix,
$$\mathbb{Q}_n^t = (Q_{n,0}, Q_{n-1,1}, \ldots, Q_{0,n}),
$$
where
$$Q_{n-k,k}(x,y)\,u = \binom{n}{k}\, \partial_x^{n-k}\,\partial_y^{k}\,
(a^{n-k}\, c^{k}\,u), \quad 0\le k \le n.$$

From the previous Lemma, $Q_{n-k,k}$ has exact degree $n$ or
$Q_{n-k,k}\equiv 0$, for $0\le k\le n$. Assume that $Q_{n-k,k}\equiv
0$, then $\partial_x^{n-k}\,\partial_y^{k}\, (a^{n-k}\, c^{k}\,u)
=0$, and therefore $a=0$ or $c=0$ which contradicts (\ref{phireg}).
Moreover,
$$\langle u, Q_{n-k,k}(x,y)\, x^{n-l}\,y^{l}\rangle =
(-1)^n \, \binom{n}{k}\, \langle a^{n-k} \, c^{k}\,u,
\partial_x^{n-k}\,\partial_y^{k}\, (x^{n-l}\, y^{l})\rangle =0,$$
if $k\neq l$. In the case $k=l$, we have
$$\langle u, Q_{n-k,k}(x,y)\, x^{n-k}\,y^{k}\rangle
\neq 0,$$ since $Q_{n-k,k}$ has exact degree $n$.

Finally, we are going to show that the polynomials $\{Q_{n-k,k},
0\le k\le n\}$ are independent modulus $\mathcal{P}_{n-1}$. Let
$\lambda_0, \lambda_1, \ldots, \lambda_n$ be constants such that
$$q(x,y) = \sum_{k=0}^n \lambda_k\, Q_{n-k,k}(x,y),$$
is a polynomial of degree less than or equal to $n-1$. Since
$$\langle u, q(x,y) p(x,y)\rangle =0,$$ for
any $p\in \mathcal{P}_{n-1}$, we have $\langle u, q(x,y)^2\rangle
=0$, and therefore $q(x,y)\equiv 0$.

Then, we get
\begin{eqnarray*}
0 &=& \langle u, q(x,y) \, x^{n-i}\, y^i\rangle = \sum_{k=0}^n
\lambda_k\,\langle Q_{n-k,k}(x,y)\, u, x^{n-i}\,y^{i}\rangle\\
&=&\lambda_i\,\langle Q_{n-i,i}(x,y)\, u, x^{n-i}\,y^{i}\rangle,
\end{eqnarray*}
 so $\lambda_i=0$, $0\le i\le n.$
\end{proof}

\bigskip

\noindent{\bf Proof of Corollary \ref{cor2}}

\begin{proof} Observe that
\begin{eqnarray*}
\langle u, \mathbb{Q}_n \mathbb{X}_n^t\rangle &=& \langle
\mathbb{Q}_n^t\, u, \mathbb{X}_n^t\rangle = \langle {\hbox{\rm
div}}^{\{n\}} (\Phi^{\{n\}}\,u), \mathbb{X}_n^t\rangle =\\
&=& (-1)^n \langle \Phi^{\{n\}}\,u,
\nabla^{\{n\}}\mathbb{X}_n^t\rangle = (-1)^n \, n!\,\langle u,
(\Phi^{\{n\}})^t\rangle,
\end{eqnarray*}
and then, the result follows.
\end{proof}

\end{document}